\documentclass[11pt]{amsart}

\usepackage{amssymb}
\usepackage{latexsym}
\usepackage{amsmath}
\usepackage{amscd}
\usepackage{array}

\def\c{{\rm C}}
\def\C{{\mathbb C}}

\def\z{{\mathbb Z}}

\def\Q{{\mathcal Q}}

\def\K{{\mathbb K}}

\begin{document}

\newtheorem{theorem}{Theorem}[section]
\newtheorem{proposition}[theorem]{Proposition}
\newtheorem{definition}[theorem]{Definition}
\newtheorem{corollary}[theorem]{Corollary}
\newtheorem{lemma}[theorem]{Lemma}
\newtheorem{conjecture}[theorem]{Conjecture}
\newtheorem{question}[theorem]{Question}

\title[algebraic quantum permutation groups]{algebraic quantum permutation groups}

\author{Julien Bichon}
\address{Laboratoire de Math\'ematiques, Universit\'e Blaise Pascal, Clermont-Ferrand II, 
Campus des C\'ezeaux, 63177 Aubi\`ere Cedex, France}
\email{Julien.Bichon@math.univ-bpclermont.fr}


\begin{abstract}
We discuss some algebraic aspects of quantum permutation groups, working over arbitrary fields. 
If $\K$ is any characteristic zero field, 
we show that there exists a universal cosemisimple Hopf algebra coacting
on the diagonal algebra $\K^n$: this is a refinement
of Wang's universality theorem for the (compact) quantum permutation group. 
We also prove a structural result for Hopf algebras having a non-ergodic coaction
on the diagonal algebra $\K^n$, on which we determine the possible group gradings
when $\K$ is algebraically closed and has characteristic zero.
\end{abstract}

\maketitle

\section{Introduction}

A remarkable fact, discovered by Wang \cite{wa}, is the existence
of a largest (universal) compact quantum group, denoted $\Q_n$, acting on the set
$[n] =\{1, \ldots ,n\}$, which is infinite if $n\geq 4$.
In view of its universal property, the quantum group $\Q_n$ is called
the quantum permutation group on $n$ points, and is seen a good  
quantum analogue of the classical permutation group $S_n$.

The compact quantum group $\Q_n$ is defined through a Hopf $\c^*$-algebra
$A_s(n)$ via an heuristic formula $A_s(n) = \c(\Q_n)$. The algebra $A_s(n)$
satisfies Woronowicz' axioms in \cite{wo} and Wang's universality theorem is stated in terms
of a coaction of $A_s(n)$ on the diagonal $\c^*$-algebra 
$\C^n$ (the algebra of functions on $[n]$). 
Canonically associated with $A_s(n)$ is 
a dense Hopf $*$-subalgebra, denoted here $A_s(n, \C)$.
It is clear from its presentation that one can define an analogous
Hopf algebra over any field $\K$, which we denote 
$A_s(n,\K)$.

In this paper we study some aspects of the Hopf algebra $A_s(n,\K)$,
which coacts on the diagonal algebra $\K^n$.
We prove that the coaction is universal in various cases (Theorem 3.8). 
In particular our main result states that in characteristic
zero $A_s(n,\K)$ is the universal cosemisimple Hopf algebra
coacting on $\K^n$, which in geometric language means that
there exists a largest linearly reductive algebraic quantum group acting on $n$ points.
This refines Wang's Theorem 3.1 in \cite{wa}.

We also prove a structural result for Hopf algebras coacting
on $\K^n$ in a non-ergodic manner. Then if $\K$ is algebraically
closed and has characteristic zero, we give a general description of
cocommutative cosemisimple quotients of $A_s(n,\K)$, which exactly correspond
to group gradings on the diagonal algebra $\K^n$.
It is worthwhile to note that  the classification problem for 
 group gradings over various classes of associative algebras (matrix algebras, triangular
algebras, incidence algebras...)
has been intensively studied in recent years: see \cite{bz,vz,ps}.

Quantum permutation groups have been studied in connection
with subfactor theory, free probability theory and the classification
problem for quantum groups.
We refer the reader to the survey paper \cite{bbc2} for an overview of 
these recent results. As well as presenting some new results, we hope
that the present paper might serve as a friendly algebraic introduction
for entry
to these developments.

The paper is organized as follows.
Section 2 consists of notations and preliminaries, with a short review of Hopf algebra
coactions and universal coactions. In Section 3 we study Hopf algebra coactions
on the diagonal algebra $\K^n$, and prove various universality results 
 for $A_s(n,\K)$. Section 4 is devoted to non-ergodic coactions.
 Group gradings on diagonal algebras are described in Section 5, leading
 to the classification of cosemisimple cocommutative quotients of $A_s(n,\K)$, if 
 $\K$ has characteristic zero and is algebraically closed.

\section{Notations and preliminaries}

\subsection{Diagonal algebras}
Throughout the paper $\K$ is field.
The diagonal algebra $\K^n$, $n \in \mathbb N^*$, is always equipped with its canonical basis
$e_1, \ldots , e_n$, with
$$e_ie_j = \delta_{ij} e_i, \quad 1 = \sum _{i=1}^n e_i$$
The algebra $\K^n$ is identified with the algebra of $\K$-valued
functions on the set $[n] = \{1, \ldots , n\}$.
Under this identification the idempotent $e_i$ is the characteristic
function of the subset $\{i\}$.   

\subsection{Hopf algebras and comodules}
We assume that the reader has some familiarity with bialgebras, Hopf 
algebras and their comodules, for which the books \cite{mo,ks}
are convenient references.
Unless otherwise indicated,
the comultiplication, counit and antipode of a Hopf algebra
are always denoted by $\Delta$, $\varepsilon$ and $S$ respectively.

Let $B$ be a bialgebra and let $x_{ij}$, $i,j \in [n]$, be some elements of $B$.
Recall that the matrix $x = (x_{ij})\in M_n(B)$ is said to be \textbf{multiplicative}
if $$\Delta(x_{ij}) = \sum_{k=1}^n x_{ik} \otimes x_{kj}, \quad \varepsilon(x_{ij}) = \delta_{ij}$$
The matrix $x$ is multiplicative if and only the linear map
$$\beta : \K^n \longrightarrow \K^n \otimes B, \quad \beta(e_i)= \sum_{k=1}^n e_k \otimes x_{ki},$$
endows $\K^n$ with a $B$-comodule structure. 

We will be interested in particular types of multiplicative matrices. 
We say that a multiplicative matrix $x =(x_{ij}) \in M_n(B)$ is \textbf{semi-magic}
if the following two families of relations hold in $B$. 
 $$ 
 \leqno(2.1) \quad x_{ki}x_{kj} = \delta_{ij} x_{ki}, \  \forall i,j,k \in [n]$$
$$ \leqno(2.2) \quad \sum_{k=1}^n x_{ik} =1, \ \forall i \in [n]$$
We say that the multiplicative matrix $x$ is \textbf{magic} if it is
semi-magic and if furthermore the following relations hold in $B$.
$$ 
 \leqno(2.3) \quad x_{ik}x_{jk} = \delta_{ij} x_{ik}, \  \forall i,j,k \in [n]$$
$$ \leqno(2.4) \quad \sum_{k=1}^n x_{ki} =1, \ \forall i \in [n]$$

\subsection{Cosemisimple Hopf algebras}
Recall that a Hopf algebra $H$ is said to be cosemisimple if its comodule
category is semisimple. This is equivalent to saying that there exists
a linear map $h : H \longrightarrow \K$ (the Haar measure)
such that 
$$h(1) = 1, \quad ({\rm id}_H \otimes h) \circ \Delta = h(-)1_H = (h \otimes {\rm id}_H) \circ \Delta$$
A CQG algebra is Hopf $*$-algebra (over $\C$) having all its finite dimensional comodules
 equivalent to unitary ones. A CQG algebra is automatically 
 cosemisimple, but there exist cosemisimple
 Hopf algebras that do  not admit any CQG algebra structure (for example
 the quantum group $SL_q(2)$ for $q$ non real and not a root of unity). CQG algebras
 are the algebras of representative functions on compact quantum groups, see \cite{ks}
 for more details.

\subsection{Coactions on algebras}
Let $A$ be an algebra and let $B$ be a bialgebra. A \textbf{coaction $B$ on the algebra 
$A$} consists of an algebra map $\beta : A \longrightarrow A \otimes B$ making
$A$ into  a $B$-comodule. One also says that $A$ is a $B$-comodule algebra.
If a linear map $\beta : A \longrightarrow A \otimes B$ endows
$A$ with a $B$-comodule structure, then it is a coaction on the algebra $A$
if and only if $A$ is an algebra in the monoidal category of $B$-comodules, that
is, the multiplication map $m: A \otimes A \longrightarrow A$ and unit map $\K \longrightarrow A$
are $B$-colinear. Coactions on algebras correspond to quantum (semi)groups actions on quantum spaces.

Let $\mathcal C$ be a subcategory of the category of bialgebras (resp. Hopf algebras).
Let $H$ be a bialgebra (resp. Hopf algebra) coacting on the algebra $A$, via
$\alpha : A \longrightarrow A \otimes H$. We say that the coaction is \textbf{universal in $\mathcal C$},
or that $H$ is \textbf{the universal bialgebra (resp. Hopf algebra) in $\mathcal C$ coacting on $A$}, if 
for any coaction $\beta : A \longrightarrow A \otimes H'$, with $H'$ being a bialgebra
(resp. Hopf algebra) in $\mathcal C$, there exists a unique bialgebra map
$f : H \longrightarrow H'$ such that     
$({\rm id}_{A} \otimes f) \circ \alpha = \beta$.

Manin first \cite{ma} proposed to construct bialgebras and Hopf algebras
by looking at universal objects (in suitable categories) coacting on well
chosen algebras (quadratic algebras in \cite{ma}).
In the case of finite-dimensional algebras, it is not difficult to show
that given a finite-dimensional algebra $A$, there exists a universal bialgebra coacting
on $A$, and hence, through Manin's Hopf envelope construction \cite{ma},
a universal Hopf algebra coacting on $A$. The problem with this Hopf algebra
is that it is not finitely generated in general, and hence  can hardly be thought of as the
algebra of functions on the quantum symmetry group of a finite quantum space.

Manin's work was continued by Wang in \cite{wa} in the framework of Woronowicz algebras,
the objects dual to compact quantum groups.
In this paper Wang showed that there exists indeed a universal
Woronowicz algebra (or equivalently a universal CQG algebra) coacting on the diagonal $\c^*$-algebra $\C^n$, denoted now $A_s(n)$.
For non-commutative $\c^*$-algebras, Wang showed that such a universal object
does not exist (as shown by the quantum $SO(3)$-groups), but studied instead
quantum symmetry groups of $\c^*$-algebras endowed with a faithful positive functional,
leading to the construction of other quite interesting quantum groups. 

In this paper we concentrate on Hopf algebra coactions on the diagonal algebra $\K^n$.
The algebra $A_s(n, \C)$, the canonically defined dense CQG
subalgebra of Wang's $A_s(n)$, is the universal $\C$-algebra generated
by the entries of a magic matrix, and it is clear that such a definition
works over any field. We are especially interested in knowing
if $A_s(n, \K)$ is still universal in an appropriate category.

\section{The universal coaction}

Let $\K$ be an arbitrary field. We study Hopf algebra coactions on the diagonal algebra
$\K^n$. We begin with bialgebra coactions, for which we
have the following basic result. The proof uses standard arguments, and is left to the reader,
who might also consult \cite{wa}.

\begin{proposition}
 Let $B$ be a bialgebra and let
 $\beta : \K^n \longrightarrow \K^n \otimes B$ be a right comodule structure 
 on $\K^n$, with
 $$\beta(e_i) = \sum_{k=1}^n e_k \otimes x_{ki}$$
 Then $\beta$ is an algebra map (and hence a coaction on the algebra $\K^n$) if and only if the
 matrix $x=(x_{ij})$ is semi-magic.
\end{proposition}

Having this proposition in hand, we see that there indeed exists a 
universal bialgebra coacting on $\K^n$, which is the universal algebra
generated by the entries of a semi-magic matrix.   

At the Hopf algebra level, Manin's Hopf envelope of the previous bialgebra furnishes
a universal Hopf algebra coacting on $\K^n$. However this Hopf algebra is certainly too
big (not finitely generated), and 
we believe that the good object is $A_s(n, \K)$, defined as follows.

\begin{definition}
The algebra $A_s(n,\K)$ is the algebra presented by generators $u_{ij}$, $i,j \in [n]$, subject to the relations making $u=(u_{ij})$ a magic matrix.
\end{definition}

Here is the first basic result regarding $A_s(n,\K)$. The proof is left to the reader.

\begin{proposition}
 The algebra $A_s(n,\K)$ has a Hopf algebra structure defined by
 $$\Delta(u_{ij}) = \sum_{k=1}^n u_{ik} \otimes u_{kj}, \quad \varepsilon(u_{ij}) = \delta_{ij}, \quad
 S(u_{ij})=u_{ji}$$
 The formula $$\alpha(e_i) = \sum_{k=1}^n e_k \otimes u_{ki}$$
 defines a coaction of $A_s(n,\K)$ on the algebra $\K^n$.
\end{proposition}

The relationship of $A_s(n,\K)$ with the symmetric group is examined
in the next proposition, where several arguments
of \cite{wa} are used. The defining relations of magic 
matrices are those of permutation matrices, but since we are in arbitrary characteristic
we cannot claim directly here that a commutative Hopf algebra is a function algebra. 

\begin{proposition}
 There exists a surjective Hopf algebra map $\pi_n : A_s(n,\K) \longrightarrow \K(S_n)$, where
 $\K(S_n)$ is the Hopf algebra of $\K$-valued functions on the symmetric group $S_n$, such that:
 \begin{enumerate}
 \item The map $\pi_n$ is an isomorphism if and only if $n\leq 3$.
 \item The map $\pi_n$ induces an isomorphism between $A_s^c(n,\K)$, the maximal commutative
 quotient of $A_s(n,\K)$, and $\K(S_n)$.
\end{enumerate}
Moreover the algebra $A_s(n,\K)$ is non commutative and infinite-dimensional if $n\geq 4$.
\end{proposition}

\begin{proof}
 The map $\pi_n$ is defined by sending $u_{ij}$ to $p_{ij}$, the function defined
 by $p_{ij}(\sigma) =\delta_{i,\sigma(j)}$.
 It is surjective because $e_\sigma$, the characteristic function of a permutation $\sigma$,
 satisfies
 $$e_\sigma = p_{\sigma(1)1} \cdots p_{\sigma(n)n}$$
 The algebra map $\pi_n$ thus induces a surjective algebra map
 $\pi_n^c:A_s^c(n,\K) \longrightarrow \K(S_n)$. Denoting the generators of $A_s^c(n,\K)$ by $x_{ij}$, we define
 a linear map $\K(S_n) \longrightarrow A_s(n,\K)$ by sending $e_\sigma$ to 
 $x_{\sigma(1)1} \cdots x_{\sigma(n)n}$. One checks easily that this a Hopf algebra map, and that 
 it is the reciprocal isomorphism of $\pi_n^c$.
 
 If $n \geq 4$, one can reproduce Wang's argument in \cite{wa}, page 201, to see that
 $A_s(n,\K)$ has an infinite dimensional quotient given by a free product
 of non-trivial algebras.
 
 It is trivial that $A_s(2, \K)$ is commutative, while some
 slightly more involved computations, left to the reader, show that $A_s(3, \K)$ is commutative as well.
 This concludes the proof.
\end{proof}

Proposition 3.4 seems to indicate that the quantum group corresponding to 
$A_s(n,\K)$ is some kind of free version of the symmetric group $S_n$. See 
\cite{bb,bc,bbc} for probabilistic and representation theoretic
meanings of freeness.

When $\K = \C$, the Hopf algebra $A_s(n,\C)$ is a CQG algebra
(and is the dense CQG algebra of Wang's algebra $A_s(n)$),
a fact that has been already mentioned in several papers. 
We reproduce the argument here, and we note that the cosemisimplicity
holds more generally in characteristic zero.

\begin{proposition}
 If $\K$ is a characteristic zero field,
 the Hopf algebra $A_s(n,\K)$ is cosemisimple .
 If $\K = \C$, then $A_s(n,\C)$ has a Hopf $*$-algebra structure
 given by $u_{ij}^*=u_{ij}$, and is a CQG algebra.
\end{proposition}

\begin{proof}
It is straightforward to check that there exists a $\K$-algebra map
$\tau : A_s(n,\K) \longrightarrow A_s(n,K)$ such that $\tau (x_{ij})=x_{ji}$. The algebra
$A_s(n,\K)$ is defined over the ordered field $\mathbb Q$, and thus $A_s(n,\K)$ is cosemisimple
by Theorem 4.7 in \cite{bi98}. When $\K = \C$, the Hopf $*$-algebra structure
is easily defined, and the generating multiplicative matrix of $A_s(n,\C)$ being
unitary, we use [Proposition 28, p.417] in \cite{ks}  to conclude that $A_s(n,\C)$ is a CQG algebra   
\end{proof}

We now wish to study when the coaction defined in Proposition 3.3 is universal.
First we need to clarify the interactions between the various relations 
defining magic matrices. This done in the following lemma.

\begin{lemma}
 Let $H$ be a Hopf algebra, and let  $x=(x_{ij}) \in M_n(H)$ be a multiplicative matrix. 
If three of the families of relations (2.1), (2.2), (2.3), (2.4) hold for the elements $x_{ij}$, then 
 the fourth family also holds, so that $x$ is magic matrix, and $S(x_{ij}) = x_{ji}$, for all $i,j$. 
\end{lemma}

\begin{proof}
 Assume that relations (2.1), (2.2) and (2.3) hold. Then 
 it is easy to see that $xx^t =I$, the identity matrix ($x^t$ is the transpose matrix). The matrix
 $x$ is invertible with inverse $S(x)$,  hence $S(x) = x^t$. Now one gets Relations 
 (2.4) by applying the antipode to Relations (2.2). The other cases are treated with similar
 arguments.  
\end{proof}

\begin{lemma}
 Let $H$ be a Hopf algebra with bijective antipode 
 coacting on the diagonal algebra $\K^n$, with coaction
 $\beta : \K^n \longrightarrow \K^n \otimes H$. Let $x_{ij} \in H$, $i,j \in[n]$, be such that 
 $$\beta(e_i) = \sum_{k=1}^n e_k \otimes x_{ki}$$
so that the matrix $x=(x_{ij})\in M_n(H)$ is semi-magic. 
For $i \in [n]$, let $u_i=\sum_{k=1}^nx_{ki}$. Then $u_i$ is invertible  for all $i \in [n]$, and
$$S(x_{ij}) = u_i^{-1}x_{ji}, \quad {\rm and} \quad 
  S^{-1}(x_{ij}) = x_{ji}u_i^{-1}, \ \forall i,j \in [n]$$
\end{lemma}

\begin{proof}
 The matrix $x$ is multiplicative because $\beta$ endows $\K^n$ with an $H$-comodule structure, and is semi-magic by Proposition 3.1. We have
 $$(x^t)x = {\rm diag}(u_1, \ldots, u_n)$$
 and since $x$ and $x^t$ are invertible with respective inverses $S(x)$ and $S^{-1}(x)^t$, 
 the elements $u_i$ are all invertible and we have
 $$S(x) = {\rm diag}(u_1^{-1}, \ldots , u_n^{-1}) x^t \quad {\rm and} \quad
 S^{-1}(x)^t = x {\rm diag}(u_{1}^{-1}, \ldots , u_{n}^{-1})$$
 This concludes the proof of the lemma. 
\end{proof}

We are now ready
to prove the main result of the section.

\begin{theorem}
 Let $H$ be a Hopf algebra coacting on the diagonal algebra $\K^n$, with coaction
 $\beta : \K^n \longrightarrow \K^n \otimes H$. Assume that one of the following conditions hold.
 \begin{enumerate}
  \item The usual integration map $\psi : \K^n \longrightarrow \K$, $e_i \longmapsto 1$,  is $H$-colinear.
  \item $S^2 = {\rm id}_H$.
  \item $K$ has characteristic zero or characteristic $p>n$ and $H$ is cosemisimple. 
 \end{enumerate}
Then there exists a unique Hopf algebra map
$f : A_s(n,\K) \longrightarrow H$ such that
$$({\rm id}_{\K^n} \otimes f) \circ \alpha = \beta$$ 
\end{theorem}

\begin{proof}
Let $x_{ij} \in H$, $i,j \in [n]$ be as in Lemma 3.7.
We already know that the matrix $(x_{ij})$ is semi-magic, and
it is clear from the construction of $A_s(n,\K)$ that we just have to prove that it is magic.
In case (1), by the $H$-colinearity of $\psi$, we see that Relations (2.4) hold in $H$, and 
by Lemma 3.6, the matrix $x$ is magic and the theorem is proved in the first case.

Assume now that $S^2 = {\rm id}_H$.
We have $S(x_{ij}) = u_i^{-1}x_{ji} =  
  S^{-1}(x_{ij}) = x_{ji}u_i^{-1}$,  $\forall i,j \in [n]$, by Lemma 3.7.
  Let $i,j,k \in [n]$ with $i \not = j$. Then, using again Lemma 3.7, we have
  $$0 = x_{kj}x_{ki} = S(x_{kj}x_{ki}) = u_k^{-1} x_{ik} x_{jk} u_k^{-1}$$
  and hence relations (2.3) hold in $H$. We conclude using
  Lemma 3.6.

We assume now that $H$ is cosemisimple. Let $h : H \longrightarrow \K$ be the Haar measure
and let $x_{ij} \in H$, $i,j \in [n]$, be as in Lemma 3.7. Put $\alpha_i = h(u_i) = h(\sum_kx_{ki})$.
The map $\varphi = (\psi \otimes h) \circ \beta$, $\K^n \longrightarrow \K$, is $H$-colinear,
with $\varphi(e_i) = \alpha_i$. Thus the bilinear form
$\omega : \K^n \otimes \K^n \longrightarrow \K$ defined
by $\omega = \varphi \circ m$, where $m$ is the multiplication of $\K^n$, is also $H$-colinear.
We have $\omega(e_i,e_j) = \delta_{ij}\alpha_i$, and if $E= {\rm diag}(\alpha_1, \ldots, \alpha_n)$,
the $H$-colinearity of $\omega$ gives 
$x^t E x = E$, and hence $x^t E = E S(x)$. Thus we have
$$\alpha_i S(x_{ij}) = \alpha_j x_{ji}, \quad \forall i,j \in [n]$$ 
Let $I = \{ i \in [n] \ | \ \alpha_i\not = 0\}$. The set $I$ is non-empty
since
$$\sum_{i=1}^n \alpha_i = \sum_{i,k=1}^nh(x_{ki}) =nh(1)= n \not =0$$
We have
$$S(x_{ij}) = \alpha_i^{-1} \alpha_j x_{ji}, \ \forall i \in I, \forall j \in [n] $$
and since $S$ is bijective, we have $x_{ij} = 0$ for $i \in I$ and $j\not \in I$.
Also for $i \not \in I$ and $j \in I$, we have, by Lemma 3.7,
$x_{ij} = u_j S(x_{ji})=0$. 

We now concentrate on the elements $x_{ij}$, $i,j \in I$.
We wish to prove that $x_0 = (x_{ij})_{i,j \in I}$ is a magic matrix.
For $i,j \in I$, we have
$$ \Delta(x_{ij}) = \sum_{k\in I} x_{ik} \otimes x_{kj}, \quad \varepsilon(x_{ij}) = \delta_{ij}$$
and hence $x_0$ is a multiplicative matrix. It is clear 
that Relations (2.1) hold for the elements $x_{ij}$, $i,j \in I$. Also, for $i \in I$
$$1 = \sum_{k \in I}x_{ik} + \sum_{k \not \in I}x_{ik} = \sum_{k \in I}x_{ik}$$
and Relations (2.2) hold for the elements $x_{ij}$, $i,j \in I$.
Let $i,j,k \in I$, with $i \not = j$. We have
$$x_{ik}x_{jk}= \alpha_k \alpha_i^{-1}S(x_{ki})\alpha_k \alpha_j^{-1} S(x_{kj})
=\alpha_k^2 \alpha_i^{-1} \alpha_j^{-1} S(x_{kj}x_{ki}) =0$$
and hence Relations (2.3) hold for the elements $x_{ij}$, $i,j \in I$.
Now by Lemma 3.6 $x_0$ is a magic matrix and in particular we have
$$\sum_{k \in I} x_{ki} =1, \ \forall i\in I$$
Hence for $i \in I$, we have
$$u_i = \sum_{k \in I} x_{ki} + \sum_{k \not \in I} x_{ki} = \sum_{k \in I} x_{ki}=1$$
and $\alpha_i=h(u_i)=1$, $\forall i \in I$. But then
$$n= \sum_{i \in I} \alpha_i = \#I$$
Thus $I = [n]$ and the proof of Theorem 3.8 is complete. 
\end{proof}

\bigskip

Combining Theorem 3.8 and Proposition 3.5, we get 
the following universality result.

\begin{theorem}
 If $\K$ has characteristic zero, the Hopf algebra $A_s(n,\K)$ is the universal
 cosemisimple Hopf algebra coacting on the algebra $\K^n$.
\end{theorem}

This result is a refinement of
Wang's universality Theorem in \cite{wa}.
Indeed it is possible to get an equivalent version of Wang's Theorem
as an immediate corollary of Theorem 3.8.

\begin{theorem}[Theorem 3.1 in \cite{wa}] 
 The CQG algebra $A_s(n, \C)$ is the universal CQG algebra coacting
 on $\C^n$. 
\end{theorem}

\begin{proof}
 Let $H$ be a CQG algebra coacting on $\C^n$: here this means that the action is moreover
 a $*$-algebra map. The elements $x_{ij}$ in the proof of Theorem 3.8 therefore satisfy $x^*_{ij}=x_{ij}$, and the Hopf algebra morphism from the proof of Theorem 3.8
 is a $*$-algebra map.
\end{proof}

We conclude the section by some remarks and questions.
First we should mention that case (1) of Theorem 3.8
is a particular case of a general machinery developed in \cite{bi00}. We have
included it here because it is immediate using the arguments
needed to prove the other cases, and also
because the invariance of classical integration
is a natural requirement (automatic in the classical case),
which further motivates the use of $A_s(n,\K)$. 

\medskip

The corepresentation theory of 
$A_s(n,\C)$ is worked out in \cite{ba0}: it is similar to the representation theory of the algebraic
group $SO(3,\C)$.
This can be generalized in characteristic zero. 
On the other hand, we do not know if $A_s(n,\K)$ is cosemisimple in  
positive characteristic $p>n$. In this case there is always a non-trivial cosemisimple Hopf
algebra coacting on $\K^n$, the function algebra on $S_n$. So we have the following question.

\begin{question}
 Assume that $\K$ has characteristic $p>n\geq 4$. Does there exist
 a universal cosemisimple Hopf algebra coacting on the algebra $\K^n$?
\end{question}

The noncommutative cosemisimple Hopf algebras constructed in \cite{bi0} show that if 
this universal Hopf algebra exists, then it is not isomorphic to $\K(S_n)$.

\medskip

It is also clear that $A_s(n,\K)$ might be defined over any ring, with functoriality properties.
This suggests that the mod $p$ reduction $A_s(n,\mathbb Z)\rightsquigarrow A_s(n,\z/p\z)$ could be
used in some contexts (especially in the context of quantum automorphism
groups of finite graphs as in \cite{bbch}), but this idea has not been
fruitful yet.

 
 


\section{Non-ergodic coactions}

In this section we study non-ergodic coactions 
on $\K^n$. First recall
that if $H$ is  a Hopf algebra coacting on an algebra $A$, the coaction is 
said to be \textbf{ergodic} if the fixed point subalgebra
$$A^{co H} = \{ a \in A \ | \ \alpha(a) = a \otimes 1\}$$
is reduced to $\K1=\K$. Also
recall that the coaction is said to be \textbf{faithful} if $H$ is generated, as an algebra, by
the space of coefficients $\{ (\psi \otimes {\rm id}_H) \circ \alpha(a), \ a \in A, \ \psi \in A^*\}$. 

From the quantum group viewpoint, ergodic coactions
correspond to transitive actions. The following result is the analogue
of the decomposition of a classical group action into disjoint orbits.

\begin{proposition}
 Let $H$ be a Hopf algebra coacting faithfully on the algebra $\K^n$, and satisfying one the 
 assumptions of Theorem 3.8. 
 Then there exists a sequence of positive integers $m_1 \geq  \ldots \geq m_k >0$, with
 $m_1 + \cdots +m_k=n$ and  $k=\dim((\K^n)^{co H})$, with ergodic coactions of $H$ on each $\K^{m_i}$, and
 with a surjective
 Hopf algebra morphism $$A_s(m_1,\K) * \cdots * A_s(m_k, \K) \longrightarrow H$$
\end{proposition}

\begin{proof}
 It is well-known that a subalgebra of a diagonal algebra is itself diagonal, and corresponds
 to a partition of the set $[n]$. Hence we have a 
 partition $[n] = X_1 \sqcup \cdots \sqcup X_k$, with $m_i = \#X_i$, $m_1 \geq \cdots \geq m_k>0$, 
 and  we have
 $$(\K^n)^{co H} = \K f_1 \oplus \cdots \oplus \K f_k$$
 with $f_i = \sum_{j \in X_i} e_j$ being a minimal coinvariant projection, $\forall i$. Moreover we have
 $$\K^n = \K^n f_1 \oplus \cdots \oplus \K^n f_k$$
 The coinvariance of the $f_i$'s ensures that the original coaction restricts
 to coactions on the diagonal algebras $\K^nf_i = \K^{m_i}$ (with unit $f_i$), and the coactions
 are ergodic by the minimality of the $f_i$'s.
 Hence we get the announced coactions, which by Theorem 3.8 produce the announced Hopf algebra map, and the surjectivity follows 
 from the faithfulness assumption. 
\end{proof}

Although this result seems to reduce the study of coactions to ergodic ones, we have 
to say that in general it is difficult to determine
the quotients of a free product. Here is a modest application to the determination
of low degree quantum permutation groups.

\begin{corollary}
 Let $H$ be a CQG algebra coacting on $\C^5$. If the coaction is not ergodic, then
 $H$ is a Hopf $*$-algebra quotient of $A_s(4,\C)$ (listed in \cite{bb3})
 or is a Hopf $*$-algebra quotient of $\C(S_3)*\C(S_2)$.
\end{corollary}

\section{Group gradings on diagonal algebras}

Let $G$ be a group and let $A$ be an algebra. Let us recall that
a $G$-grading on $A$ consists of a vector space decomposition
$$A = \bigoplus_{g \in G} A_g$$
with $A_gA_h \subset A_{gh}$ for all $g,h \in G$.
It is then easy  to see that $1 \in A_1$.
Also it is well known that a $G$-grading is exactly
the same thing as a coaction on $A$ by  $\K[G]$, the convolution group algebra of $G$.
If we have a $G$-grading, the formula $\alpha(a) = a\otimes g$, $a \in A_g$, defines
a coaction of $\K[G]$ on $A$, and the converse follows from the cosemisimplicity of 
$A$ and the fact that its simple comodules are one-dimensional.
We freely interchange the two notions.

For a non-zero element $a \in A_g$, we write $|a|=g$. 
We say that the grading is faithful if the set
$$S = \{ g \in G \ | \ \exists a \in A \ {\rm with} \ |a|=g\}$$
generates $G$ as a group. Faithful $G$-gradings correspond to faithful coactions as in the previous section.  Of course we are only interested in faithful gradings.
We say that the grading is ergodic if $A_1=\K1$, which means that 
the corresponding coaction is ergodic.

\medskip

In this section we determine the possible group gradings on diagonal algebras. The conclusion will be given
in Proposition 5.2.

We begin with a lemma, very similar to Lemma 5 in \cite{ps}.

\begin{lemma}
 Let $G$ be a group and assume that $A=\K^n$ has a faithful $G$-grading.
 \begin{enumerate}
 \item If $A_g \not = 0$, then $g$ has finite order.
 \item If the grading is ergodic, the group $G$ is abelian.
\end{enumerate}

\end{lemma}

\begin{proof}
 Let $a \in A$ with $|a|=g$. Let $l >0$ be such that 
 there exist $\lambda_0, \ldots , \lambda_{l-1} \in \K$ such that
$a^l = \sum_{i=0}^{l-1} \lambda_ia^i$. Then
 $$a^l\otimes g^l = \sum_{i=0}^{l-1} \lambda_i a^i \otimes g^i$$ 
 If $g$ does not have finite order, then $a^l=0$ and $a=0$ ($\K^n$ has no non-zero nilpotents),
 which contradicts the assumption on $a$.  
 
 We assume now that the grading is ergodic. Let $a,b \in A$
 with $|a|=g$ and $|b|=h$. The elements $g,h \in G$ have finite order
 and hence by Lemma 5 in \cite{ps} the elements $a$ and $b$ are invertible.
 Thus $ab=ba$ is a non-zero element in $A_{gh} \cap A_{hg}$, and $gh=hg$.
 The grading being faithful, we conclude that $G$ is abelian. 
\end{proof}

\begin{proposition}
Let $G$ be a group. Assume that $\K$ has characteristic zero and is algebraically closed. Then the following are equivalent.
\begin{enumerate}
 \item There exists a faithful $G$-grading on $\K^n$.
 \item There exists a family of transitive abelian groups
 $G_i \subset S_{m_i}$, $i=1, \ldots, k$, with $m_1+ \cdots + m_k=n$, and a surjective group morphism
 $$G_1 * \cdots * G_k \longrightarrow G$$
\end{enumerate}
If these conditions hold, the grading is ergodic if and only if $G\subset S_n$ is a transitive
abelian group.
\end{proposition}

\begin{proof}
 Let $G \subset S_n$ be an abelian group. We begin by constructing
 a $G$-grading on $\K^n$. The action of $G$ on $[n]$ induces
 a coaction $\K^n \longrightarrow \K^n \otimes \K(G)$. Combined with the Hopf algebra
 isomorphisms $\K(G) \backsimeq \K[\widehat{G}] \backsimeq \K[G]$, this gives a 
 $G$-grading on $\K^n$, which is ergodic if and only if $G$ is transitive.

 Assume now that condition (2) holds. The previous construction
 gives a (transitive) $G_i$-grading on $\K^{m_i}$, and hence a 
 $G_1 * \cdots * G_m$-grading on $\K^n = \K^{m_1} \oplus \cdots \oplus \K^{m_k}$.
 This gives finally a faithful $G$-grading on $\K^n$.
 
 Assume now that $\K^n$ has a faithful $G$-grading. If the grading is ergodic,
 the group $G$ is abelian by Lemma 5.1, and the $\K[G]$-coaction gives a $\K(G)$-coaction, and hence
 a transitive $G$-action on $[n]$. The grading is faithful and hence $G \subset S_n$ is a transitive
 abelian group. In general, we have, by Proposition 4.1 and its proof,
 a surjective Hopf algebra map 
 $A_s(m_1,\K) * \cdots * A_s(m_k, \K) \longrightarrow \K[G]$
 and the image of $A_s(m_i)$ is a group algebra $\K[G_i]$ that
 coacts ergodically on $\K^{m_i}$. By the ergodic case we know that each $G_i\subset S_{m_i}$
 is a transitive abelian group, and we are done.
\end{proof}

\begin{corollary} Assume that $\K$ has characteristic zero and is algebraically closed.
 Every cocommutative and cosemisimple Hopf algebra quotient of
 $A_s(n,\K)$ is isomorphic to $\K[G]$, where the group $G$ is a quotient
 of a free product of transitive abelian groups.
 
 Every cocommutative Hopf $*$-algebra quotient of
 $A_s(n,\C)$ is isomorphic to $\C[G]$, where the group $G$ is a quotient
 of a free product of transitive abelian groups.
 
\end{corollary}

\begin{proof}
 This follows directly from the previous result because a cosemisimple
 and cocommutative Hopf algebra is a group algebra, and also a  cocommutative
 CQG algebra is a group algebra.
\end{proof}

\end{document}